\documentclass[11pt,a4paper]{article}
\usepackage{amsfonts,amsmath,amstext,amssymb}
\usepackage{dsfont}
\usepackage{theorem}    %para cambiar el estilo de los teoremas
\usepackage{a4wide}

\newtheorem{lema}{Lemma}[section]

\newtheorem{teor}[lema]{\bf Theorem}

\theorembodyfont{\rmfamily}

\title{A dual rigidity of the sphere and the hyperbolic plane}

\author{Magdalena Caballero and Rafael M. Rubio \\[6mm]
 Departamento de Matemeticas, Campus de Rabanales, \\[0.5mm] Universidad de C\'ordoba, 14071 C\'ordoba, Spain,\\[0.5mm] E-mails\textup{: \texttt{magdalena.caballero@uco.es},\quad  \texttt{rmrubio@uco.es}}}

\date{}

\begin{document}

\maketitle

\begin{abstract}
There are several well-known characterizations of the sphere as a regular surface in the Euclidean space. By means of a purely synthetic technique, we get a rigidity result for the sphere without any curvature conditions, 
nor completeness or compactness. As well as a dual result for the hyperbolic plane, the spacelike sphere in the Minkowski space.
\end{abstract}



\section{Introduction}

In 1897 Hadamard proved that any compact connected regular surface with positive Gaussian curvature in the three-dimensional Euclidean space $\mathbb{E}^3$ is a topological sphere, \cite{Hadamard}. 
His result motivated the search for conditions to conclude that such a surface is necessarily a round sphere (an Euclidean sphere). 
Two answers were given by Liebmann. The first one, in 1899 \cite{Li}, proved the rigidity of the sphere, conjectured by F. Minding in 1939, 

\begin{quote}{\it If $S$ is a compact and connected regular surface in $\mathbb{E}^ 3$ with constant Gaussian curvature $K$, then $M$ is a sphere of radius $1/\sqrt{K}$.}
\end{quote} 

The second one involves the mean curvature, \cite{Li2},

\begin{quote}{\it Any compact and connected regular surface in $\mathbb{E}^3$ with positive Gaussian curvature and constant mean curvature is a sphere.}
 
\end{quote}

Shortly after, Hilbert gave a simpler proof of the first result, \cite{Hilbert}. His ideas where used by Chern in \cite{Ch} to get a more general characterization of the sphere, concerning Weingarten surfaces. 
He obtained the previous two results by Liebmann as corollaries. 

Another result on the sphere involving the Gaussian curvature, which is a direct consequence of a result by Hopf \cite{Hopf-1926}, 
asserts that \textit{it is the only complete and simply connected regular surface in $\mathbb{E}^ 3$ with positive constant Gaussian curvature.} 

The last characterization of the sphere we will mention is the Alexandrov theorem, which assures that {\it a compact and connected regular surface of constant mean curvature in $\mathbb{E}^3$ is a sphere}, 
\cite{A}.

Now consider the Minkowski space $\mathbb{L}^3$. A regular surface in this space is called spacelike if its induced metric is Riemannian. In this setting, 
the hyperbolic plane $\mathbb{H}^2$ can be realized as one connected component of the hyperboloid of two sheets, and so it can be viewed as the spacelike sphere in $\mathbb{L}^3$. 
Analogously to the (Euclidean) sphere, the hyperbolic plane can be characterized as {\it the only spacelike regular surface in $\mathbb{L}^ 3$ which is complete, 
simply connected and with negative constant Gaussian curvature}, \cite{Hopf-1926}.

In this work, we are interested in surfaces foliated by circles. R. L\'opez proved that {\it a surface in $\mathbb{E}^3$ with constant Gaussian curvature and foliated
by pieces of circles is included in a sphere, or the planes containing the circles of the foliation are parallel}, \cite{Lopez01}. In \cite{Lopez02}, 
the same author obtained the dual result in $\mathbb{L}^3$. It states that {\it a spacelike surface in $\mathbb{L}^3$ with constant Gaussian curvature and foliated
by pieces of circles must be a portion of a hyperbolic plane, unless the planes of the foliation are parallel}.

This paper is devoted to prove natural dual characterizations of the sphere in $\mathbb{E}^3$ and the hyperbolic plane in $\mathbb{L}^3$. In our results neither the Gaussian curvature nor the mean curvature appear. 
Neither completeness nor compactness hypotheses are required. We only need a hypothesis on the intersection of the surfaces by planes.  
%In comparison with the previous characterizations, our technique is purely synthetic, i.e., without the use of parametrizations nor coordinates.

More specifically, we say that a regular surface $S$ in the Euclidean space ${\mathbb{E}^3}$ satisfies the ${\cal P}$ property if for each affine plane $\Pi$ intersecting $S$, 
the set $\Pi\cap S$ is a circle (including the degenerate case with radius zero).

Analogously, we say  that a spacelike regular surface $S$ in the Minkowski space $\mathbb{L}^3$ satisfies the ${\cal P}^*$ property if for each spacelike affine plane $\Pi$ intersecting $S$, 
the set $\Pi\cap S$ is a circle (including the degenerate case with radius zero). Notice that a circle in a spacelike affine plane $\Pi$ of $\mathbb{L}^3$ 
is the locus of the points in $\Pi$ at a constant distance from a fixed point in $\Pi$, where the distance considered is the one associated to the induced metric.

We state the following rigidity results:

\begin{teor}
Let $S$ be a connected regular surface in the Euclidean space ${\mathbb{E}^3}$ satisfying the ${\cal P}$ property, then $S$ is necessarily an Euclidean sphere. 
\end{teor}

\begin{teor}
Let $S$ be a spacelike connected regular surface in the Minkowski space ${\mathbb{L}^3}$ satisfying the ${\cal P^*}$ property, then $S$ is necessarily a hyperbolic plane. 
\end{teor}

\section{The proofs}

{\bf{Euclidean case.}}

\vspace{2mm}

Let $S$ be a surface in $\mathbb{E}^3$ satisfying the ${\cal P}$ property and let $Q\in S$ be an arbitrary point. We consider the tangent plane $T_{_Q}S$ and its normal line through $Q$, ${\cal L}$. 
We take the sheaf of affine planes %${\cal A}$ 
with axis ${\cal L}$ and we denote by $\{C_i\}_{i\in I}$ the family of circles obtained when intersecting those planes with $S$.

We consider a plane $\Pi_0$ parallel to $T_{_Q}S$ such that $C=\Pi_0\cap S$ is a non degenerate circle and we denote $P=\Pi_0\cap {\cal L}$. Then $\Pi_0\cap C_i\not=\emptyset$ for all $i\in I$, 
and the intersection points of each circle $C_i$ with $\Pi_0$ are the opposite points of a chord of $C_i$ contained in $\Pi_0$ with midpoint $P$. 
Therefore, the point $P$ must be the center of $C$ and as direct consequence the circles $C_i$ have all the same radius. Thus, the sphere given by ${\bigcup}_{i\in I} C_i$ is contained in $S$. 
We finish the proof thanks to the connectedness of $S$.

\vspace{2mm}

\noindent{\bf{Lorentzian case.}}

\vspace{2mm}

Let $S$ be a surface in $\mathbb{L}^3$ satisfying the ${\cal P^*}$ property and let $\Pi_0$ be a spacelike plane such that $C=\Pi_0\cap S$ is a non-degenerate circle. 
We denote its center by $P$ and the normal line through $P$ by ${\cal L}$. 

Firstly, we prove that ${\cal L}\cap S\neq \emptyset$. We proceed by contradiction. Let us assume ${\cal L}\cap S= \emptyset$. 
Therefore, any plane parallel to $\Pi_0$ either does not intersect $S$ or it does it in a non-degenerate circle whose interior contains a point of ${\cal L}$. 
We deduce that any line parallel to ${\cal L}$ intersects $S$ at most in one point, otherwise the intersection of $S$ and the plane generated by both lines contains a non-spacelike curve. 
Thus, we have proved that $S$ is a graph over a domain of $\Pi_0$ not intersecting ${\cal L}$ and foliated by circles. Since $S$ is spacelike, it can not be asymptotic to ${\cal L}$ or any line parallel to it, 
and so $\partial S\neq \emptyset$. If $\partial S$ does not contain a point of ${\cal L}$, 
then it contains a circle. In both cases we can find spacelike planes intersecting $S$ in a non closed curve, which is a contradiction. 

We notice that $S$ must by closed, otherwise we proceed as before to arrive to a contradiction. We take a point $Q$ at which the distance from $P$ to ${\cal L}\cap S$ is attained. 

For each $A\in C$ and for each chord perpendicular to the segment $AP$, we call its midpoint  $A_m$. 
%Since ${\cal L}$ is perpendicular to $\Pi_0$, 
We can chose the chord as close to $A$ as necessary so that the plane generated by it and the segment $A_mQ$ is spacelike, we denote it by $\Pi_{A_m}$. 
We define $\varepsilon_A$ to be the supremum (in the set of all possible chords satisfying the previous property) of the distance from $A$ to $A_m$. 

If $\varepsilon=\min_C\varepsilon_A$, we take $0<\rho<\varepsilon$ and for each $A\in C$ we consider the chord with $d(A_m,A)=\rho$. Hence, all the circles $\Pi_{A_m}\cap S$ 
have the same radius, and so there is a hyperbolic cap contained in $S$ and containing $C$. 

Finally, for each point $A\in S$ there exists a spacelike plane intersecting $S$ in a non-degenerate circle containing $A$. Therefore, there exists a hyperbolic cap contained in $S$ and containing $A$. 
We finish the proof by using a connectedness argument.

\section*{Acknowledgments}
The  authors are partially supported by the Spanish MICINN Grant with FEDER funds
MTM2010-18099.

\end{document}